\documentclass{article}[11pt]

\usepackage{amssymb}
\usepackage{amsmath}
\usepackage[mathscr]{eucal}
\usepackage{epsfig}
\usepackage{fullpage}
\usepackage{graphicx,bm}

\newtheorem{Theorem}{Theorem}

\newcommand{\ed}{\ \stackrel{d}{=} \ }

\newcommand{\TT}{\mbox{${\mathcal T}$}}

\newcommand{\Zbold}{\mbox{${\mathbb Z}$}}
\newcommand{\Tbold}{\mbox{${\mathbb T}$}}

\newcommand{\TTT}{\mbox{${\mathscr{T}}$}}

\newcommand{\bE}{{\bf E}}
\newcommand{\bP}{{\bf P}}

\newcommand{\bone}{\mathbf{1}}

\title{{\bf On the Expected Total Number of Infections for Virus Spread on a Finite Network}}

\author{{\bf Antar Bandyopadhyay}\footnote{E-Mail: antar@isid.ac.in} \\  
        {\bf Farkhondeh Sajadi}\footnote{E-Mail: farkhondeh.sajadi@gmail.com}\\
        \ \\
        Theoretical Statistics and Mathematics Unit,\\
        Indian Statistical Institute, Delhi Centre,\\
        7 S. J. S. Sansanwal Marg, \\
        New Delhi 110016\\
        INDIA}

\date{}

\begin{document}

\maketitle

\begin{abstract}
In this paper we consider a simple virus infection spread model on a finite population of $n$ agents 
connected by some neighborhood structure. Given a graph $G$ on $n$ vertices, we begin with some fixed number of 
initial infected vertices. At each discrete time step, an infected vertex tries to infect its neighbors with probability 
$\beta \in (0,1)$ independently of others and then it dies out. The process continues till all infected vertices die out. 
We focus on obtaining proper lower bounds on the expected number of ever infected vertices. We obtain a simple lower bound, 
using \textit{breadth-first search} algorithm and show that for a large class of graphs which can be classified as the ones which 
locally ``look like'' a tree in sense of the \emph{local weak convergence} \cite{AlSt04}, this lower bound gives better approximation than 
some of the known approximations through matrix-method based upper bounds \cite{DrGaMa08}. \\

\noindent
\emph{{\bf AMS 2000 Subject Classifications:}} Primary: 60K35, 05C80; secondary: 60J85, 90B15 \\

\noindent
\emph{{\bf Keywords and phrases:}} Bread-first search, local weak convergence, random $r$-regular graphs, susceptible infected removed model,
                                   virus spread.

\end{abstract}

\section{Introduction}
\label{Sec:Intro}

\subsection{Background and Motivation}
\label{SubSec:Background}
Often it is observed that the normal operation of a system which is organized in a network of 
individual machines or agents is threatened by the propagation of a harmful entity through the network.
Such harmful entities are often termed as a \emph{viruses}. 
For example the Internet, as a network is threatened by the computer viruses and worms which are self-replicating pieces of
code, that propagate in a network of computers. These codes use a number of different methods to propagate, for example, 
an e-mail virus typically sends copies of itself to all addresses in the address book of the infected machine. Weaver et. al. 
\cite{WePaStCu03} gives a good survey of different techniques of propagation for computer viruses.

In this paper we use a simple \emph{susceptible infected removed (SIR)} model which was studied by Draief, Ganesh and Massouli\'{e} in 2008 \cite{DrGaMa08}. 
In this model, each susceptible agent, can be infected by its infected neighbors at a rate, proportional to their number and remains infected
till it is removed after an unit time. While it is infected, it has the potential to infect its neighbors. 
In general, removal can correspond to a quarantining of a machine from the network or patching the machine. 
In this model, it is assumed that once a node is removed, it is ``out of the network''. That is, it can no longer be susceptible or infected. 
Such a model is justified, provided the epidemic spread happens at a much faster rate than the rate of patching of the susceptible machines.

The study of mathematical models for epidemic spread has a long history in biological epidemiology and in the study of computer viruses. One of the first 
work in this area was by Kermack and Mckendrick \cite{KerMc27}, where they established the first stochastic theory for epidemic spread. They also proved 
the existence of an epidemic threshold, which determines whether the epidemic will spread or die out. As mentioned in \cite{DrGaMa08}, earlier work mainly 
focused on finding or approximating the \emph{law of large numbers} limit where the stochastic behavior was approximated by its mean behavior and hence 
mainly studied deterministic models. More recent works \cite{BaUt04, LeUt95}, have focused on stochastic nature of the models and have tried to prove 
asymptotic distribution of the
number of survivors, using a key concept called 
\emph{basic reproductive number} $R_0$, which is defined as the expected number of secondary
infective, caused by a single primary infective. This concept of basic reproductive number is well defined under the
\emph{uniform mixing} assumption, that is, when any infective can infect any susceptible equally likely, and hence the
underlying network is given by a complete graph. For a general network, where basic reproductive number may become 
vertex dependent, it is not clear how to use this concept effectively. 
As in \cite{DrGaMa08},
in this work we would like to study this model on a general network.

\subsection{Model}
\label{SubSec:Model}
We consider a closed population of $n$ agents, connected by a network structure, given by an undirected graph $G = \left(V, E\right)$ with vertex set $V$, 
containing all the agents 
and edge set $E$. A vertex can be in either of the three states, namely, \emph{susceptible (S)}, \emph{infected (I)} or \emph{removed (R)}. At the 
beginning, the
initial set of infected vertices is assumed to be non-empty and all others are susceptible. The evolution of the epidemic is described by the following discrete 
time model: 
\begin{itemize}
\item After a unit epoch of time, 
      each infected vertex instantaneously tries to infect each susceptible neighbor with       
      probability $\beta \in \left(0,1\right)$ 
      independent of all others. 
\item Each infected vertex is removed from the network after an unit time. 
\end{itemize}
Mathematically, at an integer multiple of unit time, say $t$, if a susceptible vertex $v$ has $I_v\left(t\right)$ neighbors who are infected, then the 
probability
of $v$ being infected instantaneously is $1 - \left(1 - \beta\right)^{I_v\left(t\right)}$ and each susceptible vertex gets infected independently. Also an
infected vertex remains in the network only for an unit time, after that it tries to infect its susceptible neighbors and then it is immediately removed. 

As pointed out by \cite{DrGaMa08}, this is a simple model, falling in the class of models known as Reed-Frost Models, where infection period is 
deterministic and is same for every vertex. It is worth noting that the evolution of the epidemic can be
modeled as a Markov chain. 

It is interesting to note here that, the model is essentially same as the i.i.d. Bernoulli bond percolation model with parameter $\beta$ \cite{Gri99}. 
This is because the set of ever infected (or removed) 
vertices is same as the union of connected open components of i.i.d. bond percolation on $G$, containing all the initial infected vertices. Although for 
percolation, it is customary to work
with an infinite graph $G$. If $G$ is the complete graph $K_n$, then this model is fairly well studied 
in literature and is known as the \emph{binomial random graph}, also known as 
Erd\"{o}s-R\'{e}nyi random graph \cite{Jan, Bola01}.

Like in \cite{DrGaMa08}, our goal is to study the total 
number of vertices that eventually become infected (and hence removed)
without specifying the underlying network. In \cite{DrGaMa08}, the authors derived an explicit upper bound of the expected number of vertices ever infected
which depends on both the size of the network as well as the infection rate $\beta$. This bound also needed an assumption of ``small'' value for $\beta$. 
Unfortunately, the work \cite{DrGaMa08} did not provide any 
indication whether the derived upper bound is a good approximation of
the quantity of interest. In this work we derive a simple lower bound of the expected number of vertices ever infected which works for every infection rate
$0 < \beta < 1$. Our lower bound is based on the \emph{breadth-first search (BFS)} algorithm and hence easily computable for any general finite network $G$. We also prove that, under certain
assumptions on the qualitative behavior of the underlying graph, namely if $G$ 
``\emph{locally looks like a tree}'' in the sense of Aldous and Steele \cite{AlSt04} 
\emph{local weak convergence}, then
our lower bound is asymptotically exact for ``small'' $\beta$, thus providing
a good approximation when the network is ``large''. As we will see later, 
for such graphs $G$, the range we cover for
$\beta$ always includes the range in which the upper bound obtained in \cite{DrGaMa08} holds and in all these
cases, the upper bound over estimates the expected total number of infections.

\subsection{Outline}
\label{SubSec:Outline}
In the following section, we state and prove our main results.
Section \ref{Sec:Example} gives several examples where our lower bound holds and gives asymptotically correct answer. 
Finally in Section \ref{Sec:Discussion} we summarize the merits of our work and indicate some of its limitations as well.

\section{Main Results and Proofs}
\label{Sec:Results}
We will denote by $Y^{G, I}$, the total number of vertices 
ever infected when the epidemic runs on a network $G$ and the infection starts at the vertices in
$I \subseteq V$. Note that $Y^{G, I}$ implicitly depends on the size of the network. 
In Subsection \ref{SubSec:Results-One} we present the results, when
the epidemic starts with only one infected vertex.
We generalize these results  for epidemic starting with more than one infection, which are presented in Subsection \ref{SubSec:Results-Many}. 
In both cases, our results relay on a specific search algorithm, known as \emph{breadth-first search (BFS)}. We briefly describe the algorithm 
here. 

\begin{quote}
{\tt
\begin{itemize}

\item[Step-0] Input graph $G$ with a linear ordering of its vertices say 
              $V := \left\{v_0, v_1, v_2, \cdots, v_{n-1}\right\}$. Let
              $T \leftarrow \left\{v_0\right\}$ and $N \leftarrow \left\{ v_0 \right\}$.
              
\item[Step-1] Write $N = \left\{v_{i_1}, v_{i_2}, \cdots, v_{i_r} \right\}$ for some $r \geq 1$
              such that $i_1 < i_2 < \cdots < i_r$.

\item[Step-2] For $l=1$ to $r$ 
              find all neighbors $u$ of $v_{i_l}$ which are not in $T$, put
              \[
              N' \leftarrow N' \cup \left\{u \,\Big\vert\, u \sim v_{i_l} \text{\ and\ } u \not\in T \,\right\} 
              \]
              and update $T$ as
              \[
              T \leftarrow T \cup  \left\{u \,\Big\vert\, u \sim v_{i_l} \text{\ and\ } u \not\in T \,\right\} \,.
              \] 

\item[Step-3] Update $N \leftarrow N'$.
              
\item[Step-4] Go to Step-1 unless vertex set of $T$ is
              same as that of $V$.
              
\item[Step-5] Stop with output $T$ as the BFS spanning tree with root $v_0$. 

\end{itemize}
}
\end{quote}

Note that the BFS spanning tree is not necessarily unique, 
it depends on the starting point $v_0$ which is typically
called the root and also it depends on the ordering of the vertices in which the exploration of neighbors is done
in {\tt Step-2}. Also note that iff $G$ is a tree to start with then, BFS spanning tree is just
itself. 


\subsection{Starting with Only One Infected Vertex}
\label{SubSec:Results-One}
Our first result gives a lower bound of the expected total number of vertices ever infected starting with exactly one infected vertex.
\begin{Theorem}
\label{Thm:LB}
Let $G$ be an arbitrary finite graph and $v_{0} \in V$ be a fixed vertex of it. 
Let $T$ be a spanning tree of the connected component of $G$ containing the vertex $v_0$
and rooted at $v_0$. Let $Y^{T, \{v_0\}}$ be  
the total number of vertices 
ever infected when the epidemic runs only on $T$ and starting with exactly one infection at $v_0$. 
Then   
\begin{equation}
\bE\left[Y^{T, \left\{v_0\right\}}\right] \leq \bE\left[Y^{G, \left\{v_0\right\}}\right] \,\,\, \text{for all} \,\,\, 0 < \beta < 1 \,.
\label{Equ:LB-Tree}
\end{equation}
Moreover, if  $\TT$ is a BFS spanning tree of the connected component of $v_0$ rooted at $v_0$, then 	
\begin{equation}
\bE\left[Y^{T, \left\{v_0\right\}}\right] \leq \bE\left[Y^{{\mathcal T}, \left\{v_0\right\}}\right] \leq 
\bE\left[Y^{G, \left\{v_0\right\}}\right] \,\,\, \text{for all} \,\,\, 0 < \beta < 1 \,.
\label{Equ:BFS-Tree}
\end{equation}
\end{Theorem}

\begin{proof}
%
Suppose $G = \left(V, E\right)$ where $V$ is the set of vertices and $E$ is the set of edges and let $H = \left(V, E'\right)$ where $E' \subseteq E$. 
So $H \subseteq G$, is a spanning sub-graph of $G$. Note that $v_0$ is a vertex in both $H$ and $G$.
Let $\left(X_e\right)_{e \in E}$ be i.i.d. $\mbox{Bernoulli}\left(\beta\right)$ random variables indexed by the edges of the graph $G$. 
We consider the random graphs $G_{\beta} := \left(V_{\beta}, E_{\beta}\right)$ and $H_{\beta} := \left( V_{\beta}, E'_{\beta}\right)$
with the same vertex set $V_{\beta} = V$ and
the random sets of edges $E_{\beta} := \left\{e \in E \,\big\vert\, X_e = 1 \,\right\}$ and
$E'_{\beta} := \left\{e \in E' \,\big\vert\, X_e = 1 \,\right\}$. Note that $H_{\beta}$ is a spanning sub-graph of $G_{\beta}$. 
Let $C^{G, v_0}$ and $C^{H, v_0}$ be the connected components of the vertex $v_0$ in $G_{\beta}$ and $H_{\beta}$ respectively. From definition
$C^{H, v_0} \subseteq C^{G, v_0}$.

Now it follows from the definition of the infection spread model that $\left\vert C^{G, v_0} \right\vert \ed Y^{G, \left\{v_0\right\}}$ and
$\left\vert C^{H, v_0} \right\vert \ed Y^{H, \left\{v_0\right\}}$.
So to prove equation (\ref{Equ:LB-Tree}) observe that
\[
     \bE\left[Y^{T, \left\{v_0\right\}}\right]
=    \bE\left[ \left\vert C^{T, \left\{v_0\right\}} \right\vert \right]
\leq \bE\left[ \left\vert C^{G, \left\{v_0\right\}} \right\vert \right]
=    \bE\left[Y^{G, \left\{v_0\right\}}\right] \,.
\]

For the second part, we note that if $T$ is a spanning tree of $G$ with root $v_0$, then 
$d_G\left(v, v_0\right) \leq $ $d_T\left(v, v_0\right)$
for all $v \in V$, where $d_G$ and $d_T$ are the graph distance functions on $G$ and $T$ respectively. Moreover, the BFS algorithm preserves the distances, 
so if $\TT$ is a BFS spanning tree with root
$\left\{v_0\right\}$ then we must have
\[d_G\left(v, v_0\right) = d_{{\mathcal T}}\left(v, v_0\right)\] for
all $v \in V$. Thus $d_{{\mathcal T}}\left(v, v_0\right) \leq d_T\left(v, v_0\right)$ for all $v \in V$. 
Now from the model description, it follows that for any spanning tree $T$ with root $v_0$ we have
\[
\bE\left[Y^{T, \left\{v_0\right\}}\right] = \sum_{v \in V} \beta^{d_T\left(v, v_0\right)} \,.
\]
So we conclude that
\[
\bE\left[Y^{T, \left\{v_0\right\}}\right] = \sum_{v \in V} \beta^{d_T\left(v, v_0\right)} \leq \sum_{v \in V} 
\beta^{d_{{\mathcal T}}\left(v, v_0\right)} = \bE\left[Y^{{\mathcal T}, \left\{v_0\right\}}\right] \,,
\]
as $0 < \beta < 1$.
\end{proof}

Let $\text{LB}^{G, \left\{v_0\right\}} := \bE\left[Y^{{\mathcal T}, \left\{v_0\right\}}\right]$ be the lower bound 
obtained through BFS algorithm for a BFS spanning tree $\TT$ of $G$, rooted at $v_0$. Then  
from the proof of Theorem \ref{Thm:LB} we get that
\begin{equation}
\text{LB}^{G, \left\{v_0\right\}} = \sum_{v \in V} \beta^{d_{G}\left(v, v_0\right)} \,,
\label{Equ:LB-Universal}
\end{equation}
which is free of the choice of  the BFS spanning tree. Later, we will see that, this helps us to
generalize the lower bound for epidemic starting with more than one infected vertex.
We also note that $\text{LB}^{G, \left\{v_0\right\}}$ can be easily computed using the
breadth-first search algorithm described earlier.

Our next result shows that if we have a ``large'' finite graph $G$ on $n$ vertices and the epidemic starts with exactly
one infected vertex $v_0$, such that any cycle containing $v_0$ is ``relatively large'', that is of order $\Omega\left(\log n\right)$, 
then the lower bound $\text{LB}^{G, \left\{v_0\right\}}$ 
given above, is asymptotically same as the exact quantity 
$\bE\left[ Y^{G, \left\{ v_0 \right\}} \right]$. 

To state the result rigorously, we use the following graph theoretic notations. 
Given a graph $G$, a fixed vertex $v_0$ of $G$ and $d \geq 1$,  let $V_d\left(G\right)$ be the set of vertices of $G$ which are at 
a \emph{graph distance} at most $d$ from $v_0$ in $G$. Let $N_d\left(G, v_0\right)$ be the induced sub-graph of $G$ on the
vertices $V_d\left(G\right)$.  

\begin{Theorem}
\label{Thm:Large-Local-Girth}
Let $\left\{(G_n,v_0^n)\right\}_{n \geq 1}$ be a sequence of rooted connected graphs on $n$-vertices with roots $\left\{v_0^n\right\}_{n \geq 1}$ 
such that there exists a sequence $\alpha_n = \Omega\left(\log n\right)$ with $N_{\alpha_n}\left(G_n, v_0^n\right)$ is a tree for all $n \geq 1$. 
Then, there exists $0 < \beta_0 \leq 1$, such that for all $0 < \beta < \beta_{0}$
\begin{equation}
\frac{\bE\left[Y^{G_n, \left\{v_0^n\right\}}\right]}{\text{LB}^{G_n, \left\{v_0^n\right\}}} 
\longrightarrow 1 
\,\,\, \text{as} \,\,\, n \rightarrow \infty \,.
\label{Equ:LB-Exact-for-Large-Local-Girth}
\end{equation}
\end{Theorem}

\begin{proof}
Let $\TT_n$ be a BFS spanning tree rooted at $v_0^n$ of the graph $G_n$ and as defined earlier and let
$\text{LB}^{G_n, \left\{v_0^n\right\}} = \bE\left[Y^{{\mathcal T}_n, \left\{v_0^n\right\}}\right]$. 
Then
\begin{eqnarray}
\text{LB}^{G_n, \left\{v_0^n\right\}}
    & \leq & \bE\left[Y^{G_n, \left\{v_0^n\right\}}\right] \nonumber \\
    & \leq & \bE\left[Y^{N_{\alpha_n}\left(G_n, v_0^n\right), \left\{v_0^n\right\}}\right] 
             + \bE\left[Y^{N_{\alpha_n}\left(G_n, v_0^n\right), \left\{v_0^n\right\}}\right]
             \times \beta^{\alpha_n} \times n \nonumber \\
    & \leq & \text{LB}^{G_n, \left\{v_0^n\right\}} + \text{LB}^{G_n, \left\{v_0^n\right\}} \times \beta^{\alpha_n} \times n \,,\label{Equ:Bound-1}
\end{eqnarray}
Note that the first term  of the second inequality in \eqref{Equ:Bound-1}  is the expected number of infected nodes within an $\alpha_{n}$ neighbourhood of the initial 
infective $v_0^n$. The second term there is an upper bound of the expected number of nodes which may become infected by these neighbourhood infectives. 
Consider the nodes which are on the boundary of  $\alpha_{n}$ neighbourhood of $v_0^n$, that is the infected vertices in $G_{n}$ after $\alpha_{n}$ units of 
time starting with one infected at vertex $v_{0}^{n}$.  Since we have assumed that $N_{\alpha_n}\left(G_n, v_0^n\right)$ is a tree, so
these nodes have probability $\beta^{\alpha_{n}}$ to get infected after $\alpha_n$ units of time. 
But the number of nodes outside the neighbourhood is bounded by $n-\bE\left[Y^{N_{\alpha_n}\left(G_n, v_0^n\right), \left\{v_0^n\right\}}\right] \leq n$. 
Therefore an upper bound for the expected number of nodes which may become
infected by these neighbourhood infectives is
\[\bE\left[Y^{N_{\alpha_n}\left(G_n, v_0^n\right), \left\{v_0^n\right\}}\right]
             \times \beta^{\alpha_n} \times n \,.\]
Also the last inequality follows from the fact that $N_{\alpha_n}\left(G_n, v_0^n\right)$ is a tree and hence is a subtree of $\TT_n$. 
This proves (\ref{Equ:LB-Exact-for-Large-Local-Girth}) since by assumption $\alpha_n = \Omega\left(\log n\right)$. 
\end{proof}

Although the assumption in the above theorem, may seem to be very restrictive, it is satisfied in many examples including the $n$-cycle 
(see Subsection \ref{SubSec:Cycle}).
The method of the proof on the other hand, helps us generalize the result for a large class of graphs including certain random graphs. 

Following Aldous and Steele \cite{AlSt04}, we say a sequence of rooted random or deterministic graphs $\left\{(G_n,v_0^n)\right\}_{n \geq 1}$
with roots $\left\{v_0^n\right\}_{n \geq 1}$ converges to a random or deterministic graph $\left(G_{\infty}, v_0^{\infty}\right)$ in the sense of 
\emph{local weak convergence (l.w.c)} and write
$\left(G_n, v_0^n\right) \xrightarrow{l.w.c.} \left(G_{\infty}, v_0^{\infty}\right)$
if for any $d \geq 1$,
\begin{equation}
\bP\left( N_d\left(G_n, v_0^n\right) \cong N_d\left(G_{\infty}, v_0^{\infty}\right) \right) 
\longrightarrow 1 \,\,\, \text{as} \,\, n \rightarrow \infty \,.
\label{Equ:LWC}
\end{equation}
Note that for a sequence deterministic graphs, (\ref{Equ:LWC}) means that the event occurs for ``large"' enough $n$. 

\begin{Theorem}
\label{Thm:Limit-for-Bounded-Degree-Graphs} 
Let $\left\{(G_n,v_0^n)\right\}_{n \geq 1}$ be a sequence of rooted connected deterministic or random graphs with deterministic
or randomly chosen roots $\left\{v_0^n\right\}_{n \geq 1}$. Suppose that for each $G_n$ the maximum degree of a vertex is bounded by 
a fixed constant, namely $\Delta$. 
Suppose there is a rooted deterministic or random tree $\TTT$ with root $\phi$ such that
\begin{equation}
\left(G_n, v_0^n\right) \xrightarrow{l.w.c.} \left(\TTT, \phi\right) \,\,\, \text{as} \,\,\, n \rightarrow \infty \,.
\label{Equ:LWC-Tree}
\end{equation}
Let $\text{LB}^{G_n, \left\{v_0^n\right\}} := \bE\left[Y^{\TT_n, \left\{v_0^n\right\}}\right]$ where $\TT_n$ is a BFS spanning tree rooted at 
$v_0^n$ of the graph $G_n$. \\
Then for $\beta < \frac{1}{\Delta}$
\begin{equation}
\left( \bE\left[Y^{G_n, \left\{v_0^n\right\}}\right] - \text{LB}^{G_n, \left\{v_0^n\right\}} \right) 
\longrightarrow 0 \,\,\, \text{as} \,\,\, n \rightarrow \infty \,. 
\label{Equ:Limit-Exact}
\end{equation}
Moreover for $\beta < \frac{1}{\Delta}$ we have
\begin{equation} 
  \lim_{n \rightarrow \infty} \text{LB}^{G_n, \left\{v_0^n\right\}}
= \lim_{n \rightarrow \infty} \bE\left[Y^{G_n, \left\{v_0^n\right\}}\right] 
= \bE\left[Y^{\mathscr{T}, \phi}\right] \,.
\label{Equ:Limit-Answer}
\end{equation}
\end{Theorem}

\begin{proof}
Let $\TT_n$ be a BFS spanning tree rooted at $v_0^n$ of the graph $G_n$ and also as defined earlier let
$\, \,\text{LB}^{G_n, \left\{v_0^n\right\}} = \bE\left[Y^{{\mathcal T}_n, \left\{v_0^n\right\}}\right]$. 
Fix $d \geq 1$ and $E_n$ be the event $\left[ N_d\left(G_n, v_0^n\right) \cong N_d\left(\TTT, \rho \right) \right]$. Therefore from Theorem \ref{Thm:LB}
\begin{equation}
\text{LB}^{G_n, \left\{v_0^n\right\}} \leq \bE\left[Y^{G_n, \left\{v_0^n\right\}}\right] = \bE\left[Y^{G_n, \left\{v_0^n\right\}} \bone_{E_n} \right] +                                                                   \bE\left[Y^{G_n, \left\{v_0^n\right\}} \bone_{E_n^c} \right] \,.
\label{Equ:Basic-UB-LB}
\end{equation}
Now under our assumption, the degree of any vertex of $G_n$ is bounded by $\Delta$ and $\beta < \frac{1}{\Delta}$, 
so using Theorem 2.3 of \cite{DrGaMa08} we have
\begin{equation}
\bE\left[Y^{G_n, \left\{v_0^n\right\}} \bone_{E_n^c} \right] \leq \frac{1}{1 - \beta \Delta} \bP\left(E_n^c\right) \,.
\label{Equ:UB-Second-Part}
\end{equation}
Further note that if $E_n$ occurs, $N_d\left(G_n, v_0^n\right)$ is a tree rooted at $v_0^n$ and thus on $E_n$, 
$N_d\left(G_n, v_0^n\right)$ is a sub-tree of $\TT_n$. So
\[      Y^{N_d\left({\mathcal T}_n, v_0^n\right), \left\{v_0^n\right\}} \bone_{E_n} 
\leq Y^{{\mathcal T}_n, \left\{v_0^n\right\}} \bone_{E_n} \,.\] 
Hence we have
\begin{eqnarray}
\bE\left[Y^{G_n, \left\{v_0^n\right\}} \bone_{E_n} \right] 
 & \leq &  \bE\left[ Y^{N_d\left({\mathcal T}_n, v_0^n\right), \left\{v_0^n\right\}} \bone_{E_n} \right]
           + \beta^d \bE\left[ Y^{G_n, \partial_d^* N_d\left(G_n, v_0^n\right)} \right] \nonumber \\
 & \leq &  \bE\left[  Y^{N_d\left({\mathcal T}_n, v_0^n\right), \left\{v_0^n\right\}} \bone_{E_n} \right]
           + \beta^d \frac{1}{1 - \beta \Delta} \bE\left[  \left\vert \partial_d^* N_d\left(G_n, v_0^n\right) \right\vert \right] \nonumber \\
 & \leq &  \text{LB}^{G_n, \left\{v_0^n\right\}} + \beta^d \frac{1}{1 - \beta \Delta} \bE\left[Y^{G_n, \left\{ v_0^n \right\}}\right] \nonumber \\
 & \leq &  \text{LB}^{G_n, \left\{v_0^n\right\}} + \beta^d \frac{1}{\left( 1 - \beta \Delta \right)^2} \label{Equ:UB-First-Part} \,,
\end{eqnarray}
where $\partial_d^* N_d\left(G_n, v_0^n \right)$ denotes the infected vertices in $G_n$ after $d$ units of time starting with one
infected vertex $v_0^n$. For the first inequality, 
note that on the event $E_n$ we have $N_d\left(G_n, v_0^n\right)$ is a tree and thus on $E_n$ each vertex in $\partial_d^* N_d\left(G_n, v_0^n \right)$
has exactly $\beta^d$ probability to get infected after $d$ units of time starting with one infected vertex at $v_0^n$.
In the second and the last inequalities, we use Theorem 2.3 of \cite{DrGaMa08}.

So finally combining (\ref{Equ:Basic-UB-LB}), (\ref{Equ:UB-First-Part}) and (\ref{Equ:UB-Second-Part}) we get that for $\beta < \frac{1}{\Delta}$
and for any $d \geq 1$ we have
\begin{equation} 
\left( \bE\left[Y^{G_n, \left\{v_0^n\right\}}\right] - \text{LB}^{G_n, \left\{v_0^n\right\}} \right)
\leq 
\beta^d \frac{1}{\left( 1 - \beta \Delta \right)^2} + \frac{1}{1 - \beta \Delta} \bP\left(E_n^c\right) \,.
\label{Equ:Difference-Exact-to-LB}
\end{equation}
Now under assumption (\ref{Equ:LWC-Tree}), we have $\mathop{\lim}\limits_{n \rightarrow \infty} \bP\left(E_n^c\right) = 0$ so we conclude that for any 
$d \geq 1$
\begin{equation}
\limsup_{n \rightarrow \infty} \left( \bE\left[Y^{G_n, \left\{v_0^n\right\}}\right] - \text{LB}^{G_n, \left\{v_0^n\right\}} \right) 
\leq
\beta^d \frac{1}{\left( 1 - \beta \Delta \right)^2} \,.
\end{equation}
This proves (\ref{Equ:Limit-Exact}) by taking $d \rightarrow \infty$ as $\beta < \frac{1}{\Delta}$.

Now for proving (\ref{Equ:Limit-Answer}), we first observe that from (\ref{Equ:LWC-Tree}) the degree of any vertex of $\mathscr{T}$ is also bounded by $\Delta$. So 
using Theorem 2.3 of \cite{DrGaMa08} we have for $\beta < \frac{1}{\Delta}$
\[
\bE\left[ Y^{N_d\left(\mathscr{T}, \rho\right), \left\{ \rho \right\}} \right] \leq \frac{1}{1 - \beta \Delta} \,.
\]
Moreover from the definition, 
$Y^{N_d\left(\mathscr{T}, \rho\right), \left\{ \rho \right\}} \uparrow 
Y^{\mathscr{T}, \left\{\rho\right\}}$ as $d \rightarrow \infty$. So
by the Monotone Convergence Theorem we have 
\begin{equation}
\lim_{d \rightarrow \infty} \bE\left[ Y^{N_d\left(\mathscr{T}, \rho\right), \left\{ \rho \right\}} \right] 
= \bE\left[ Y^{\TTT, \left\{\rho\right\}} \right] \leq \frac{1}{1 - \beta \Delta} < \infty \,.
\label{Equ:UB-on-Limit}
\end{equation}
Thus for fixed $\epsilon > 0$ we can find $d \geq 1$ such that
\begin{equation}
\left\vert \bE\left[ Y^{\TTT, \left\{\rho\right\}} \right] - \bE\left[ Y^{N_d\left(\mathscr{T}, \rho\right), \left\{ \rho \right\}} \right] \right\vert < \epsilon
\label{Equ:Limit-d-eps-1}
\end{equation}
and
\begin{equation}
\beta^d \frac{1}{\left(1-\beta \Delta\right)^2} < \epsilon \,.
\label{Equ:Limit-d-eps-2} 
\end{equation}
The last inequality holds as $\beta < \frac{1}{\Delta} < 1$. 
Further, as degree of any vertex of $\TTT$ is bounded by $\Delta$ so arguing similar to the derivation of the equation (\ref{Equ:UB-Second-Part}) we conclude
\begin{equation}
  \bE\left[ Y^{N_d\left(\mathscr{T}, \rho\right), \left\{ \rho \right\}} \right] 
- \bE\left[ Y^{N_d\left(\mathscr{T}, \rho\right), \left\{ \rho \right\}} \bone_{E_n} \right]
= \bE\left[ Y^{N_d\left(\mathscr{T}, \rho\right), \left\{ \rho \right\}} \bone_{E_n^c} \right]
\leq \frac{1}{1 - \beta \Delta} \bP\left(E_n^c\right) \,.
\label{Equ:TTT-d-to-whole}
\end{equation}
Also, arguing similar to the derivation of the equation (\ref{Equ:Difference-Exact-to-LB}) we get
\begin{eqnarray} 
\left\vert \bE\left[Y^{G_n, \left\{v_0^n\right\}} \right] -  
\bE\left[ Y^{N_d\left(G_n, v_0^n\right), \left\{ v_0^n \right\}}  \bone_{E_n} \right] \right\vert 
& \leq & 
\beta^d \frac{1}{\left( 1 - \beta \Delta \right)^2} + \frac{1}{1 - \beta \Delta} \bP\left(E_n^c\right) \nonumber \\
& \leq & \epsilon + \frac{1}{1 - \beta \Delta} \bP\left(E_n^c\right)\,,
\label{Equ:Difference-Exact-to-Approx}
\end{eqnarray}
where the last equality follows from (\ref{Equ:Limit-d-eps-2}).
Finally, 
\begin{eqnarray*}
           \left\vert \bE\left[Y^{G_n, \left\{v_0^n\right\}} \right] -   \bE\left[ Y^{\TTT, \left\{\rho\right\}} \right] \right\vert 
& \leq &   \left\vert \bE\left[Y^{G_n, \left\{v_0^n\right\}} \right] -   \bE\left[ Y^{N_d\left(G_n, v_0^n\right), \left\{ v_0^n \right\}}  \bone_{E_n} \right] \right\vert \\
&      & + \left\vert \bE\left[ Y^{N_d\left(G_n, v_0^n\right), \left\{ v_0^n \right\}}  \bone_{E_n} \right] - \bE\left[ Y^{N_d\left(\mathscr{T}, \rho\right), \left\{ \rho \right\}} \right] \right\vert \\
&      & + \left\vert \bE\left[ Y^{N_d\left(\mathscr{T}, \rho\right), \left\{ \rho \right\}} \right] - \bE\left[ Y^{\TTT, \left\{\rho\right\}} \right] \right\vert \\
& \leq & 2 \epsilon + \frac{2}{1 - \beta \Delta} \bP\left(E_n^c\right) \,,
\end{eqnarray*}
where the last inequality follows from the equations (\ref{Equ:Limit-d-eps-1}), (\ref{Equ:Limit-d-eps-2}), 
(\ref{Equ:TTT-d-to-whole}) and (\ref{Equ:Difference-Exact-to-Approx}) and also observing the fact that
$\bE\left[ Y^{N_d\left(G_n, v_0^n\right), \left\{ v_0^n \right\}}  \bone_{E_n} \right] = 
\bE\left[ Y^{N_d\left(\mathscr{T}, \rho\right), \left\{ \rho \right\}} \bone_{E_n} \right]$.
Now under our assumption (\ref{Equ:LWC-Tree}) we have $\bP\left(E_n\right) \longrightarrow 1$. So 
we conclude that
\begin{equation} 
  \lim_{n \rightarrow \infty} \bE\left[Y^{G_n, \left\{v_0^n\right\}}\right] 
= \bE\left[Y^{\TTT, \left\{\rho\right\}}\right] \,.
\label{Equ:Limit-Answer-1}
\end{equation}

Thus using (\ref{Equ:Limit-Exact}), it follows that
\[
  \lim_{n \rightarrow \infty} \text{LB}^{G_n, \left\{v_0^n\right\}}
= \lim_{n \rightarrow \infty} \bE\left[Y^{G_n, \left\{v_0^n\right\}}\right] 
= \bE\left[Y^{\mathscr{T}, \left\{\rho\right\}}\right] \,.
\]
This completes the proof. 
\end{proof}

An immediate and interesting application of the above theorem is the following result 
which gives an explicit formula for the limit of epidemic spread on a randomly 
selected $r$-regular graph when the infection starts from an
randomly chosen vertex.  
\begin{Theorem}
\label{Thm:r-Reg-Graphs}
Suppose $G_n$ is a graph selected uniformly at random from the set of all $r$-regular graphs on $n$ vertices where we assume
$n r$ is an even number. Let $v_0^n$ be an uniformly selected vertex of $G_n$. Then for $\beta < \frac{1}{r}$
\begin{equation}
\lim_{n \rightarrow \infty} \bE\left[ Y^{G_n, \left\{v_0^n\right\}} \right] = \frac{1 + \beta}{1 - \left(r-1\right) \beta} \,.
\label{Equ:r-Reg-Ans}
\end{equation}
\end{Theorem}
We note that in this case, the upper bound given in \cite{DrGaMa08} is $\frac{1}{1 - r \beta}$ when $\beta < \frac{1}{r}$ which is
strictly bigger than the exact answer given in (\ref{Equ:r-Reg-Ans}).

\begin{proof}
It is known \cite{Jan, AlSt04}
that if $G_n$ is a graph selected uniformly at random from the set of all $r$-regular graphs on $n$ vertices, where $n r$ is even and $v_0^n$ be a 
randomly selected vertex of $G_n$ then
\begin{equation}
\left(G_n, v_0^n\right) \xrightarrow{l.w.c.} \left(\Tbold_r, \phi\right) \,,
\label{Equ:LWC-r-Reg-Graphs}
\end{equation}
where $\Tbold_r$ is the infinite $r$-regular tree with root say $\phi$. The result then follows 
from Theorem \ref{Thm:Limit-for-Bounded-Degree-Graphs} and equation (\ref{Equ:r-Reg-Tree-Ans}). \end{proof}

\subsection{Starting with More than One Infected Vertex}
\label{SubSec:Results-Many}
Now suppose instead of one infected vertex, we start with $k$ infected vertices given by 
$I := \left\{v_{0,1}, v_{0,2}, \cdots, v_{0,k} \right\}$. The following theorem gives a
lower bound similar to that of Theorem \ref{Thm:LB}. 
\begin{Theorem}
\label{Thm:LB-Many} 
Let $G$ be an arbitrary finite graph and $I := \left\{v_{0,j}\right\}_{j=1}^k$ be a fixed set of $k$ vertices. 
Let $T$ be a spanning forest of the connected components of $G$ containing the vertices in $I$ with
exactly $k$ trees which are rooted at the vertices in $I$. Then   
\begin{equation}
\bE\left[Y^{T, I}\right] \leq \bE\left[Y^{G, I}\right] \,\,\, \text{for all} \,\,\, 0 < \beta < 1 \,.
\label{Equ:LB-Forest}
\end{equation}
Moreover, if  $\TT$ is a \emph{breath-first-search spanning forest} of the connected components of $G$ 
containing the vertices in $I$ with
exactly $k$ trees which are rooted at the vertices in $I$ then 	
\begin{equation}
\bE\left[Y^{T, I}\right] \leq \bE\left[Y^{{\mathcal T}, I}\right] \leq \bE\left[Y^{G, I}\right] \,\,\, \text{for all} \,\,\, 0 < \beta < 1 \,.
\label{Equ:BFS-Forest}
\end{equation}
\end{Theorem}
Given a finite labeled graph $G$ and a fixed set of vertices 
$I = \left\{v_{0,j}\right\}_{j=1}^k$ 
of it, by a \emph{breath-first-search spanning forest} 
of the connected components of $G$ 
containing the vertices in $I$
with exactly $k$ trees which are rooted at the vertices in $I$,
we mean a spanning forest of $G$ with exactly $k$ connected components which are rooted at the vertices 
$\left\{v_{0,1}, v_{0,2}, \cdots, v_{0,k}\right\}$, 
that are obtained through the \emph{breath-first-search} algorithm, starting at some vertex $v \in I$ and assuming that all the vertices
$\left\{v_{0,1}, v_{0,2}, \cdots, v_{0,k}\right\}$ are at the same level. Alternately, we can consider
a new graph $G^*$ which is same as $G$ except it has one
``artificial'' vertex, say $v^*$ which is connected to the vertices 
$v_{0,1}, v_{0,2}, \cdots, v_{0,k}$ through $k$ 
``artificial'' edges and we perform the BFS algorithm on $G^*$ 
starting with the vertex $v^*$, to obtain a 
BFS spanning tree, say $\TT^*$ of $G^*$ rooted at $v^*$. Then a \emph{breath-first-search spanning forest} of $G$ 
with exactly $k$ trees which are rooted at the vertices $\left\{v_{0,1}, v_{0,2}, \cdots, v_{0,k}\right\}$ is given by the
forest $\TT^* \setminus \left\{ v^* \right\}$. This alternate description, is quite useful in practice. Note that if 
$\left\{ \TT_i \right\}_{1 \leq i \leq k}$ are the $k$ connected components, rooted respectively at $\left\{ v_{0,1}, v_{0,2}, \cdots, v_{0,k} \right\}$ 
of $\TT$, a breath-first-search spanning forest of the connected components of $G$ containing the vertices in $I$, then the following identity holds for 
every $\beta \in \left(0,1\right):$
\begin{equation}
\bE\left[Y^{{\mathcal T}, I}\right] = \sum_{i=1}^k  \bE\left[Y^{{\mathcal T}_i, I}\right] = 
\frac{\bE\left[Y^{{\mathcal T}^*, \left\{v^*\right\}}\right] - 1}{\beta} \,.
\label{Equ:One-to-Many-Fundamental}
\end{equation}
Using the above identity, we can now generalize all the results of the previous section for 
epidemic spread starting with more than one infected vertex. 

We write $\text{LB}^{G, I}$ for $\bE\left[Y^{{\mathcal T}, I}\right]$ which is the lower bound 
of $\bE\left[Y^{G, I}\right]$
for starting with $k$ infected vertices given by $I$. 
Observe
that from equation (\ref{Equ:One-to-Many-Fundamental}) we can write
\begin{equation}
\text{LB}^{G, I} = \sum_{i=1}^k  \bE\left[Y^{{\mathcal T}_i, I}\right]\,,
\label{Equ:One-to-Many-Representation}
\end{equation}
where $\TT = \mathop{\cup}\limits_{i=1}^k \TT_i$ is as above. It is worth nothing here that the lower bound $\text{LB}^{G,I}$ does not depend on the 
choice of $\TT$ but the representation given in 
equation (\ref{Equ:One-to-Many-Representation}) uses a specific choice of $\TT$.

\begin{Theorem}
\label{Thm:Large-Local-Girth-Many}
Let $\left\{(G_n, I_n)\right\}_{n \geq 1}$ be a sequence of graphs where each $G_n$ has $k$-roots  
given by $I_n := \left\{v_{0,1}^n, v_{0,2}^n, \cdots, v_{0,k}^n \right\}$ 
such that there exists a sequence $\alpha_n = \Omega\left(\log n\right)$ with 
$N_{\alpha_n}\left(G_n, I_n\right) := \mathop{\cup}\limits_{j=1}^k N_{\alpha_n}\left(G_n, v_{0,j}^n\right)$ is a forest with $k$ components. 
Then there exists $0 < \beta_0 \leq 1$, such that for all $0 < \beta < \beta_{0}$
\begin{equation}
\frac{\bE\left[Y^{G_n, I_n} \right]}{\text{LB}^{G_n,I_n}} \longrightarrow 1 
\,\,\, \text{as} \,\,\, n \rightarrow \infty \,.
\label{Equ:LB-Exact-for-Large-Local-Girth-Many}
\end{equation}
\end{Theorem}
The proof of this result is similar to that of 
Theorem \ref{Thm:Large-Local-Girth} and follows from the
identity (\ref{Equ:One-to-Many-Fundamental}). The details are thus omitted. 

Our next result is parallel to the Theorem \ref{Thm:Limit-for-Bounded-Degree-Graphs} which needs 
a generalization of the concept of local weak convergence which was introduced by 
W\"{a}stlund \cite{Wa11}. 

We will say a sequence of random or deterministic graphs $\left\{G_n\right\}_{n \geq 1}$
with $k$ roots given by $I_n := \left\{v_{0,1}^n, v_{0,2}^n, \cdots, v_{0,k}^n\right\}$, $n \geq 1$ converges to a random or deterministic graph 
$G_{\infty}$ with $k$-roots 
say $I_{\infty} := \left\{v_{0,1}^{\infty}, v_{0,2}^{\infty}, \cdots, v_{0,k}^{\infty}\right\}$
in the sense of 
\emph{local weak convergence (l.w.c)} and write
$\left(G_n, I_n\right) \xrightarrow{l.w.c.} \left(G_{\infty}, I_{\infty}\right)$
if for any $d \geq 1$
\begin{equation}
\bP\left( N_d\left(G_n, v_{0,j}^n\right) \cong N_d\left(G_{\infty}, v_{0,j}^{\infty}\right) 
\,\,\, \text{for all} \,\,\, 1 \leq j \leq k\right) \longrightarrow 1 \,\,\, \text{as} \,\, n \rightarrow \infty \,.
\label{Equ:LWC-Many}
\end{equation}
Note that for a  sequence deterministic graphs, (\ref{Equ:LWC-Many}) means that the event occurs for ``large"' enough $n$.

\begin{Theorem}
\label{Thm:Limit-for-Bounded-Degree-Graphs-Many} 
Let $\left(G_n\right)_{n \geq 1}$ be a sequence of deterministic or random graphs. Suppose 
each $G_n$ has deterministic
or randomly chosen $k$ roots given by $I_n := \left\{v_{0,1}^n, v_{0,2}^n, \cdots, v_{0,k}^n \right\}$ and maximum degree of each $G_n$ is 
bounded by $\Delta$. 
Suppose $\TTT := \mathop{\cup}\limits_{j=1}^k \TTT_j$ is a forest with $k$ rooted tress with roots
$I_{\infty} := \left\{ \phi_1, \phi_2, \cdots, \phi_k \right\}$. We assume that
\begin{equation}
\left(G_n, I_n\right) \xrightarrow{l.w.c.} \left(\TTT, I_{\infty}\right) \,\,\, \text{as} \,\,\, n \rightarrow \infty \,.
\label{Equ:LWC-Tree-Many}
\end{equation} 
Then for $\beta < \frac{1}{\Delta}$
\begin{equation}
\left( \bE\left[Y^{G_n, I_n}\right] - \text{LB}^{G_n,I_n} \right) \longrightarrow 0 \,, 
\label{Equ:Limit-Exact-Many}
\end{equation}
as $n \rightarrow \infty$.
Moreover
\begin{equation} 
  \lim_{n \rightarrow \infty} \text{LB}^{G_n,I_n}
= \lim_{n \rightarrow \infty} \bE\left[Y^{G_n, I_n}\right] 
= \bE\left[Y^{\mathscr{T}, I_{\infty}}\right]
= \sum_{j=1}^k \bE\left[ Y^{\mathscr{T}_j, \left\{\phi_j\right\}} \right] \,.
\label{Equ:Limit-Answer-Many}
\end{equation}
\end{Theorem}

\begin{proof}
For each $n \geq 1$ as done above we define a new rooted graph $G_n^*$ with artificial vertex 
$v_n^*$ which is connected to the the $k$-roots in $I_n$ of $G_n$ through $k$ artificial edges. 
Also we consider $\TTT^*$ defined similarly with an artificial root $\phi^*$ connecting to
$\left\{ \phi_1, \phi_2, \cdots, \phi_k \right\}$. Then our assumption of local weak convergence
(\ref{Equ:LWC-Tree-Many}) is equivalent to
\begin{equation}
\left( G_n^*, v_n^* \right) \xrightarrow{l.w.c.} \left(\TTT^*, \phi^*\right) \,.
\end{equation}
This together with the relation (\ref{Equ:One-to-Many-Fundamental}) and Theorem 
\ref{Thm:Limit-for-Bounded-Degree-Graphs} 
completes the proof. \end{proof}

It is worth noting that in case $\left\{\TTT_j\right\}_{1 \leq j \leq k}$ are i.i.d. (if they are
random) or isomorphic (if they are constant) then equation (\ref{Equ:Limit-Answer-Many}) can be
reformulated as
\begin{equation} 
  \lim_{n \rightarrow \infty} \text{LB}^{G_n,I_n}
= \lim_{n \rightarrow \infty} \bE\left[Y^{G_n, I_n}\right] 
= \bE\left[Y^{\mathscr{T}, I_{\infty}}\right]
= k \, \bE\left[ Y^{\mathscr{T}_1, \left\{\phi_1\right\}} \right] \,.
\label{Equ:Limit-Answer-Many-2}
\end{equation}

As in the case of starting with one infected vertex, the following theorem is an immediate application of the above results.
\begin{Theorem}
\label{Thm:r-Reg-Graphs-Many}
Suppose $G_n$ is a graph selected uniformly at random from the set of all $r$-regular graphs on $n$ vertices where we assume
$n r$ is an even number. Let $I_n := \left\{v_{0,1}^n, v_{0,2}^n, \cdots, v_{0,k}^n \right\}$ be $k$ uniformly and independently 
selected vertices of $G_n$. Then for $\beta < \frac{1}{r}$
\begin{equation}
\lim_{n \rightarrow \infty} \bE\left[ Y^{G_n, I_n} \right] 
= k \frac{1 + \beta}{1 - \left(r-1\right) \beta} \,.
\label{Equ:r-Reg-Ans-Many}
\end{equation}
\end{Theorem}

\begin{proof}
Since the vertices in $I_n$ are selected unformly at random so from \cite{AlSt04} we have
\begin{equation}
\left(G_n, I_n\right) \xrightarrow{l.w.c.} \left( \TTT_r, I_{\infty}\right) \,,
\label{Equ:LWC-r-Reg-Graphs-Many}
\end{equation}
where $I_{\infty} := \left\{ \phi_1, \phi_2, \cdots, \phi_k \right\}$ and
$\TTT_r$ is a forest with $k$ infinite $r$-regular tree with roots in $I_{\infty}$.
The result then follows 
from Theorems \ref{Thm:Limit-for-Bounded-Degree-Graphs-Many} and \ref{Thm:r-Reg-Graphs}. 
\end{proof}

Once again we note that in this case, the upper bound $\frac{k}{1 - r \beta}$ given in
\cite{DrGaMa08} for $\beta < \frac{1}{r}$, is
strictly bigger than the exact answer given in (\ref{Equ:r-Reg-Ans-Many}) and the gap
increases with $k$, the initial number of infections.

\section{Examples}
\label{Sec:Example}

\subsection{Tree}
\label{SubSec:Tree}
If $G$ is a tree and the epidemic starts with only one infected vertex say $\phi$ which we call the root, then from the
construction of the lower bound it is clear that 
$\text{LB}^{G, \left\{\phi\right\}} = \bE\left[Y^{G, \left\{\phi\right\}}\right]$. In certain cases
one can find explicit formula for this quantity. Two such examples are discussed below.

\paragraph{Regular Tree}
Consider a rooted $r$-array tree ($r \geq 2$), with height $m$, denote it by $T\left(r,m\right)$. 
In $T\left(r,m\right)$ every internal vertex except the root $\phi$ has degree $r$. A vertex $v$ is said to be an internal vertex if it has a neighbor 
which is not on the unique path from $v$ to $\phi$. We assume that the 
degree of the root $\phi$ is $\left(r-1\right)$. 
Let $\mu_{m}:=\bE[Y^{T\left(r,m\right), \left\{\phi\right\}}]$. Note that 
the total number of vertices in $T\left(r, m\right)$ is $\frac{\left(r-1\right)^{m+1}-1}{r-2}$. 
Now, to calculate the exact value of $\mu_m$ we note that
\begin{equation}
\mu_{m}=1 + \left( r - 1 \right) \beta \mu_{m-1}
\end{equation}
which gives the formula
\begin{equation}
\mu_{m} = \frac{\left[\left( r - 1\right) \beta\right]^{m+1} - 1}{\left( r - 1\right) \beta - 1} \,.
\label{r-tree}
\end{equation}
As $T\left(r,m\right)$ is a tree so the lower bound is exact, that is, 
$\text{LB}^{T\left(r,m\right), \left\{\phi\right\}} = \mu_m$. 
Now the upper bound from \cite{DrGaMa08} is $\frac{1}{1 - r \beta}$ for $\beta <\frac{1}{r}$.
If $\beta < \frac{1}{r}$ then 
by Theorem \ref{Thm:Limit-for-Bounded-Degree-Graphs} we get
\begin{equation}
\bE\left[ Y^{T\left(r\right), \left\{\phi\right\}} \right] = 
\lim_{m \rightarrow \infty} \mu_m = \frac{1}{1 - \left( r - 1 \right) \beta} \,,
\label{limmu}
\end{equation}
where $T\left(r\right)$ is the rooted infinite $r$-regular tree, where each vertex except the root $\phi$
has degree $r$ and the degree of the root is $\left(r - 1 \right)$. 

We observe a gap between the lower bound which in this case agrees with $\mu_m$ to that
of the upper bound obtained in \cite{DrGaMa08}.

Now let $\Tbold_r$ be the infinite $r$-regular tree where each vertex including the root has
degree $r$. Such a tree can be viewed as disjoint union of $r$ rooted infinite $r$-regular trees
whose roots are joint to the root, say $\phi$ of $\Tbold_r$. Thus
from (\ref{limmu}) we get that for $\beta < \frac{1}{r}$
\begin{equation}
 \text{LB}^{{\mathbb T}_r, \left\{\phi\right\}} 
= \bE\left[Y^{{\mathbb T}_r, \left\{\phi\right\}} \right]
= 1 + \frac{r \beta}{1 - \left( r - 1 \right) \beta}
= \frac{1 + \beta}{1 - \left( r - 1 \right) \beta} \,.
\label{Equ:r-Reg-Tree-Ans}
\end{equation}

\paragraph{Galton-Watson Tree}
Consider a Galton-Watson branching process starting with one individual. Let the mean of the offspring distribution be $c > 0$. We denote the random 
tree generated by this process as $\text{GW}\left(c\right)$ with root $\phi$. Once
again, as discussed above since $\text{GW}\left(c\right)$ is a tree, so
$\text{LB}^{\text{GW}\left(c\right), \left\{\phi\right\}} 
= \bE\left[Y^{\text{GW}\left(c\right), \left\{\phi\right\}}\right]$. Now in this case, the 
epidemic process starting with only one infection at $\phi$, is a Galton-Watson branching process starting with one individual as the root and 
with mean of the new progeny distribution being $\beta c$. So in particular if $\beta < \frac{1}{c}$ then
from standard branching process theory 
$\bE\left[Y^{\text{GW}\left(c\right), \left\{\phi\right\}}\right] < \infty$ and
equals $\frac{1}{1 - \beta c}$ \cite{AthNey04}.

\subsection{Cycle}
\label{SubSec:Cycle}
Cycle graph is a graph that consists of a single cycle. We denote the cycle with $n$ vertices by $C_{n}$. For simplicity we assume $n$ is odd and then
from the BFS algorithm, it is immediate that starting with one infected individual, say 
at $v_0^n$, we have
\begin{equation}
\text{LB}^{C_n, \left\{v_0\right\}} = 1 + 2 \left( \beta + \beta^2 + \cdots +\beta^{\frac{n-1}{2}}\right)
\end{equation}
which converges to $\frac{1+\beta}{1-\beta}$ as $n \rightarrow \infty$ for any $0 < \beta < 1$. Now
it is clear from the definition that 
\begin{equation}
\left(C_n, v_0^n\right) \xrightarrow{l.w.c.} \left(\Zbold, 0\right) \,.
\end{equation} 
Thus using Theorem \ref{Thm:Limit-for-Bounded-Degree-Graphs} we conclude that if $\beta < \frac{1}{2}$
then 
\begin{equation} 
  \lim_{n \rightarrow \infty} \text{LB}^{C_n, \left\{v_0^n\right\}}
= \lim_{n \rightarrow \infty} \bE\left[Y^{C_n, \left\{v_0^n\right\}}\right] 
= \frac{1+\beta}{1-\beta} \,.
\label{Equ:Limit-Answer-Cycle}
\end{equation}
In fact this holds for any $0 < \beta < 1$. This is because for a cycle graph, the assumption in Theorem \ref{Thm:Large-Local-Girth} holds for $\alpha_n = n/3$. 
Thus from the proof of Theorem \ref{Thm:Large-Local-Girth} we conclude that the 
equation (\ref{Equ:Limit-Answer-Cycle}) holds for any $0 < \beta < 1$.

Now if the epidemic starts with $k$ initial infected vertices given by 
$I_n := \left\{ v_{0,1}^n, v_{0,2}^n, \cdots, v_{0,k}^n\right\}$ which are uniformly distributed, then it is easy to see that 
\begin{equation}
\left(C_n, I_n\right) \xrightarrow{l.w.c.} \left(\Zbold_j, 0\right)_{1 \leq j \leq k} \,,
\label{Equ:Cycle-LWC-Many}
\end{equation} 
where $\Zbold_j$ is just a copy of $\Zbold$. Then by Theorem \ref{Thm:Limit-for-Bounded-Degree-Graphs-Many} we conclude that for $0 < \beta < \frac{1}{2}$,
\begin{equation} 
  \lim_{n \rightarrow \infty} \text{LB}^{C_n,I_n}
= \lim_{n \rightarrow \infty} \bE\left[Y^{C_n, I_n}\right] 
= k \, \frac{1+\beta}{1-\beta} \,.
\label{Equ:Limit-Answer-Many-Cycle}
\end{equation}
As earlier we can use Theorem \ref{Thm:Large-Local-Girth-Many} with $\alpha_n = O\left(n\right)$ to conclude that 
(\ref{Equ:Limit-Answer-Many-Cycle}) holds for all all $0 < \beta < 1$.

\subsection{Generalized Cycle}
\label{SubSec:GCycle}
Suppose in a cycle graph we choose randomly without replacement $2 m$ vertices and connect these vertices by joining edges between them where $m \geq 1$ is 
fixed. We call this graph a Generalized Cycle and denote it by $\text{GC}\left(n,m\right)$. 
Now consider the epidemic model on this graph with one initial infected site $v_0^n$. For large enough $n$, the probability of 
having at least one of the $m$ pairs inside a neighborhood of $v_0^n$ of radius $r$ is given by
\[
1-\left(1-\frac{2r (2r+1)}{n(n-1)}\right)^{m}
\] 
which tends to zero as $n \rightarrow \infty$.
Therefore, a fixed neighborhood of the root is a tree with high probability, in fact it is isomorphic to a neighborhood of integer line. Hence by Theorem \ref{Thm:Limit-for-Bounded-Degree-Graphs} it follows that for $\beta < \frac{1}{2}$ 
\begin{equation} 
  \lim_{n \rightarrow \infty} \text{LB}^{\text{GC}\left(n,m\right), \left\{v_0^n\right\}}
= \lim_{n \rightarrow \infty} \bE\left[Y^{\text{GC}\left(n,m\right), \left\{v_0^n\right\}}\right] 
= \frac{1+\beta}{1-\beta} \,.
\label{Equ:Limit-Answer-GCycle}
\end{equation}
Similarly if we start with $k$ initial infected sites, say $I_n := \left\{ v_{0,j}^n\right\}_{j=1}^k$
which are chosen uniformly at random, then it is easy to see that
\begin{equation}
\left(\text{GC}\left(n,m\right), I_n\right) 
\xrightarrow{l.w.c.} \left(\Zbold_j, 0\right)_{1 \leq j \leq k} \,,
\label{Equ:GCycle-LWC-Many}
\end{equation}
where $\Zbold_j$ is just a copy of $\Zbold$. Thus by Theorem \ref{Thm:Limit-for-Bounded-Degree-Graphs-Many} we get
\begin{equation} 
  \lim_{n \rightarrow \infty} \text{LB}^{\text{GC}\left(n,m\right), I_n}
= \lim_{n \rightarrow \infty} \bE\left[Y^{\text{GC}\left(n,m\right), I_n}\right] 
= k \frac{1+\beta}{1-\beta} \,,
\label{Equ:Limit-Answer-Many-GCycle}
\end{equation}
when $\beta < \frac{1}{3}$, because the maximum degree in $\text{GC}\left(n,m\right)$ is $3$.

\subsection{Cube graph}
The cube graph is the graph obtained from the vertices and edges of the $3$-dimensional unit cube. We denote it by $Q_3$.  
Suppose initially only the vertex $(0,0,0)$ is infected. 
Consider a BFS spanning tree $\TT$ of $Q_3$ rooted at $(0,0,0)$. Since $Q_3$ has only $8$ vertices so
$Y^{{\mathcal T}, \left\{(0,0,0)\right\}}$ takes values $\left\{0,1,2,3,4,5,6,7\right\}$ and
\begin{eqnarray*} 
\text{LB}^{{\mathcal T}, \left\{(0,0,0)\right\}}
& = &  \bE\left[Y^{{\mathcal T}, \left\{(0,0,0)\right\}}\right] \\
& = & 1 + 3 \beta + 3\beta^{2} + \beta^{3} \\
& = & \left(1 + \beta \right)^3 \,.
\end{eqnarray*}
In general, the $d$-dimensional cube graph say $Q_d$ is a $d$-regular graph which has $n=2^d$ vertices.
Following a similar calculation as done above, one can show that for an epidemic starting at one vertex, the lower bound obtained in Theorem \ref{Thm:LB}
for the expected total number of vertices ever infected is given by $(1+\beta)^{d}$. 

In this example computation of the exact value of $\bE\left[Y^{Q_d, \left\{ (0,0,0) \right\}}\right]$ is difficult, but we note that there is a 
gap between the upper bound obtained in \cite{DrGaMa08}, namely $\frac{1}{1 - d \beta}$ which is valid only when $\beta < \frac{1}{d}$ and our lower 
bound. However this is an example which does not fall under any of the theorem we discuss in this paper and hence we are not sure if
the lower bound gives a good approximation.

\section {Discussion}
\label{Sec:Discussion}
The goal of this study has been to get a better idea of the expected total number of vertices ever infected with as little assumption as possible on 
the underlying graph $G$. Our approach has been to find an appropriate lower bound of this expectation. Although from a practical point of view, approximation
from above with an upper bound is a more conservative method. As shown in the examples given in Section \ref{Sec:Example}, the only known upper
bounds obtained in \cite{DrGaMa08} often over estimate the exact quantity. Moreover the upper bounds in
\cite{DrGaMa08} hold only for ``small'' values of the parameter $\beta$. 
For an arbitrary finite network, we have obtained a lower bound of the expectation of the number of vertices ever infected for any value of the parameter $\beta$ 
which is computable through the breadth-first search algorithm. Theorems 
\ref{Thm:Large-Local-Girth}, \ref{Thm:Limit-for-Bounded-Degree-Graphs},
\ref{Thm:Large-Local-Girth-Many} and \ref{Thm:Limit-for-Bounded-Degree-Graphs-Many} 
show that this lower bound is asymptotically exact for a large class of graphs when
$\beta$ value is ``small'', which always includes the values of $\beta$ for which the upper bounds in 
\cite{DrGaMa08} are defined. 

However, 
we would also like to mention here that even though the lower bound we present, works for any infection parameter $0< \beta < 1$, if the underlying 
graph has many loops, such as the complete graph $K_n$, then it does not necessarily give a good approximation.
To see this, consider the complete graph $K_n$ and suppose that the epidemic starts at a fixed vertex $v_0$. Then the lower bound $\text{LB}^{K_n, \left\{v_0\right\}} = 1 + \left(n-1\right) \beta$. Now, let 
$X_1$ be the number of infected vertices at time $t=1$. In this case it is easy to see that 
$X_1 \sim \text{Binomial}\left( n - 1,\beta \right)$. Let $u$ be one of $n-1-X_1$ vertices which are not infected at time $t=1$. Since $K_n$ is the complete graph, so the conditional probability of $u$ becomes infected at time $t=2$ given $X_1$ is 
$1 - \left( 1 - \beta \right)^{X_1}$. Hence
\begin{eqnarray*}
\bE\left[Y^{K_{n}, \left\{v_0\right\}}\right]
                   & \geq & 1+\left(n-1\right)\beta + \bE\left[\left(n-1-X_{1}\right) \left(1-\left(1-\beta\right)^{X_{1}}\right) \right] \\
                   &  =   & 1+\left(n-1\right)\beta + \left(n-1\right)-
                            \left(n-1\right)\left(1-\beta^2\right)^{n-1} \\
                   &      & \qquad \qquad - \left(n-1\right)\beta 
                            +\left(n-1\right)\beta \left(1-\beta\right)\left(1-\beta^2\right)^{n-2}
\end{eqnarray*}
Therefore we get 
\begin{equation}
\limsup_{n \rightarrow \infty}
\frac{\bE\left[Y^{K_{n}, \left\{v_0\right\}}\right] - \text{LB}^{K_n, \left\{v_0\right\}}}
     {\text{LB}^{K_n, \left\{v_0\right\}}} \geq \frac{1-\beta}{\beta}\,.
\end{equation}
where $\text{LB}^{K_n, \left\{v_0\right\}} :=\bE[Y^{{\mathcal T}_{n}, \left\{v_0^n\right\}}]$.

Here, it is worth mentioning that for the complete graph if we start with one infected vertex, then as discussed in Section \ref{Sec:Intro} the set of vertices 
ever infected is no other than 
an Erd\"{o}s-R\'{e}nyi random graph with parameter $n$ and $\beta$. Thus asymptotic behavior of $\bE\left[Y^{K_n, \left\{v_0\right\}}\right]$ is well understood in
the literature \cite{Jan, Bola01}.

\end{document}